\documentclass[11pt,amstex,amscd,epsfig]{article}
\usepackage{amssymb}

\def\Bbb{\mathbb}	

\def\al{\alpha}

\def\om{\omega}

\def\text#1{{\em #1}}

\def\be{\begin{equation}}
\def\ee{\end{equation}}
\def\bear{\begin{eqnarray}}
\def\eear{\end{eqnarray}}
\def\best{\begin{eqnarray*}}
\def\eest{\end{eqnarray*}}
\def\pf{{\bf Proof}: }

\renewcommand{\theequation}{\arabic{section}.\arabic{equation}}

\newtheorem{theorem}{Theorem}[section]
\newtheorem{prop}[theorem]{Proposition}

\newtheorem{lemma}[theorem]{Lemma}

\newtheorem{defn}[theorem]{Definition}

\newtheorem{remark}[theorem]{Remark}
\newenvironment{rem}{\begin{remark}\rm}{\end{remark}}

\newtheorem{example}[theorem]{Example}

\def\mspa{\vphantom{\sum}}
\def\bspa{\vphantom{\int}}

\def\non{\noindent}
\def\pf{\non {\bf Proof. }}
\def\qed{\nopagebreak \hskip .1in { $\Box$ }\penalty10000 %
\hskip\parfillskip \par  }

\def\ra{\rightarrow}

\def\r#1{\right#1}
\def\l#1{\left#1}

\def\ma#1{\mathop {#1} \limits}

\def\De{\Delta}

\def\P{{ \Bbb P}}
\def\Q{{ \Bbb Q}}

\def\wt#1{\widetilde{#1}}
\def\wh#1{\widehat{#1}}
\def\ov#1{\overline{#1}}

\def\M{{\cal M}}

\def\longra{\longrightarrow}
\def\H1{{\cal H}^1}

\def\De{\Delta}

\def\cc{{\cal C}}
\def\cg{{\cal C}_g^d}

\def\O{{\cal O}}

\def\phi{\varphi}
\def\longra{\longrightarrow}

\title{\bf Relations in the tautological ring of  $\M_g$}

\author{ Eleny-Nicoleta Ionel\thanks{partially supported by
the N.S.F. and a Sloan research fellowship} \\ University of Wisconsin \\
 Madison, WI 53706 }

\addtocounter{section}{0}
\begin{document}

\maketitle

\vskip.15in

\begin{abstract} Using a simple geometric argument, we obtain an
infinite family of nontrivial relations in the tautological ring of
$\M_g$ (and in fact that of $\M_{g,2}$). One immediate consequence
 of these relations is that the
classes $\kappa_1,\dots,\kappa_{[g/3]}$ generate the tautological ring
of $\M_g$, which has been conjectured by Faber in \cite{fa}, and
recently proven at the level of {\em cohomology} by Morita in
\cite{mo}. 
\end{abstract}

\bigskip

Throughout this paper we will assume $g\ge 2$. Let $\M_g$ be the
moduli space of nonsingular genus $g$ curves,  $p:\cc_g\ra
\M_g$  its universal curve and let $\psi\in A^1(\M_g)$ be the first chern
class of the relative dualizing sheaf $\om_\pi$. Let
$\kappa_a=\pi_*(\psi^{a+1})\in A^a(\M_g)$.  The tautological ring
$R^*(\M_g)\subset A^*(\M_g)$ is the sub-ring of the Chow ring (taken
with $\Q$ coefficients) generated by the $\kappa$ classes.

\medskip

Mumford proved in \cite{mu} that the classes $\kappa_1, \dots,
\kappa_{g-2}$ generate the tautological ring $R^*(\M_g)$. The purpose of this
paper is to show how to use similar  methods to obtain 
other relations in the tautological ring $R^*(\M_g)$. 
As a consequence, Theorem \ref{T.Fab} implies 
that the classes $\kappa_1,\dots,\kappa_{[g/3]}$ generate the
tautological ring $R^*(\M_g)$, a fact conjectured by Faber in 
\cite{fa}.

Morita \cite{mo} recently used completely different methods to obtain 
polynomial relations among $\kappa$ classes in {\em cohomology}, which
also prove Faber's conjecture in cohomology. Morita's 
methods are based on symplectic representation theory and don't seem
to extend to the Chow ring. In particular, he uses a ``crushing'' map
to induce  relations from a higher genus to a lower genus. Such crushing
map is very natural in cohomology, but it does not have a
correspondent in the Chow ring.

By contrast, our Theorem \ref{T.rel.y} uses a simple idea
to obtain relations in the tautological ring $R^*(\M_g)$.  Let
$w_{d,g}\in A^*(\M_g)$ be the class corresponding to the moduli space
of curves which can be written as degree $d$ covers of $\P^1$. Then
for each $b\ge 0$, $\kappa_{b-1} w_{d,g}$ can be computed by first
adding a marked point to the domain and then marking all the other
$d-1$ points in the same fiber, while $w_{d,g}$ can be also computed
by first marking one of the $r$ ramification points of the cover and
then marking the remaining $d-2$ points in the same fiber. Porteous
formula expresses each of these cycles in terms of
tautological classes on $\M_g$, giving 
a relation in the Chow ring for each degree $d\ge 2$ and each $b\ge
0$.

The geometric relations from Theorem \ref{T.rel.y} have a simple
generating function expression. There is no reason, a priori, why all
these relations would be given by one generating function, and
especially one in the form of an exponential of a {\em linear}
combination of $\kappa$ classes. But in Theorem \ref{T.rel.simple} we
prove that that there are some coefficients $c_{k,j}$ such that when
$b=0$ the geometric relation in $R^{g+1-2d}(\M_g)$ simply becomes
\best \l[ \;\exp\l(- \sum_{a=1}^\infty x^{a}\kappa_{a} \sum_{j=0}^{a}
c_{a,j} u^j \r) \;\r]_{x^{g+1-2d}u^d}=0 \eest Here, and throughout
this paper, the notation $ [f(t)]_{t^k} $ means the coefficient of
$t^k$ in a formal power series $f(t)$.

Furthermore,  for $b\ge 0$ the geometric relations become
\bear
\l[ \;\exp \l(-\sum_{a=1}^\infty  x^{a}\kappa_{a} \sum_{j=0}^a 
c_{a,j} u^j \r)\cdot
 \l(\kappa_{b-1}-2\ma\sum_{a=0}^\infty \kappa_{a+b} x^{a+1} 
\sum_{j=0}^{a} q_{a,j}u^{j+1} \r)\;\r]_{x^{g+2-2d}u^d}
=0
\label{eq.m.q}
\eear
where the coefficients $q_{k,j}$ are positive integers related to the
coefficients $c_{k,j}$ by the formula (\ref{q=c}). 
Theorem \ref{T.Fab} (part of Faber's Conjecture) follows then as an
easy consequence.

\medskip
The paper is organized as follows. We begin the first section  by
proving the geometric relations mentioned above (Theorem
\ref{T.rel.y}). We next 
state the main Theorem \ref{T.rel.simple} and discuss its
consequences including Faber's Conjecture (Theorem \ref{T.Fab}).  The
proof of Theorem \ref{T.rel.simple} occupies the remaining  sections of the
paper. First, in Section 2 we express the geometric relations from Theorem
\ref{T.rel.y} in terms of a universal generating function $G$ which is
the solution of a simple ODE (Theorem \ref{T.rel}).  Then 
in Section 3 we study the properties of the generating function $G$,
including a certain power series expansion for both $G$ and its
derivative. These are used in Theorem \ref{T.rel.simple} to prove that
the relations of Theorem \ref{T.rel.y} simplify to the form
(\ref{eq.m.q}). We conclude Section 3 with the proof of
 Proposition \ref{P.gen.c.aa}
which gives another simple generating function for some of the 
above relations.
\smallskip

In fact, as explained in the Appendix, the relations discussed in this
paper come from relations in the tautological ring of $\cc_g^2$
(which in turn are restrictions of relations in the Chow ring
$A^*(\ov\M_{g,2})$). Such relations are obtained by following the
approach of \cite{i}. Consider the moduli space of degree $d$, genus
$g$ covers of $\P^1$ which have a fixed ramification pattern over four 
fixed points in $\P^1$. Just as in \cite{i}, one can split the 
target $\P^1$ into two rational components, each having two of the 
four fixed points. There are several ways to split, depending which 
points are together on the same bubble. When the ramification patterns
are appropriately chosen, one obtains then relations in the Chow ring
of the Deligne-Mumford moduli space. 

In particular, the relation (\ref{eq.m.q}) is obtained from the relation
in $A^*(\cc_g)$ 
\bear\label{main.rel.c}
\l[ \;\exp \l(-\sum_{a=1}^\infty  x^{a}\kappa_{a} \sum_{j=0}^a 
c_{a,j} u^j \r)\cdot
 \l(1-2\ma\sum_{a=0}^\infty \psi_1^{a+1} x^{a+1} 
\sum_{j=0}^{a} q_{a,j}u^{j+1} \r)\;\r]_{x^{g+2-2d}u^d}
=0
\eear
after multiplying by $\psi_1^b$ and pushing forward to $\M_g$. Note
that the $\kappa$ classes in (\ref{main.rel.c}) are those pulled back 
from $\M_g$. Relation (\ref{main.rel.c}) gives a nontrivial polynomial 
expression in  $\psi_1$  and $\kappa$ classes in degree 
$a=g+2-2d$ for each $d\le (g+2)/3$. The coefficient of $\psi_1^a$ is
$-2q_{a-1,d-1}$, so nonzero. As far as we know, this is the first time
this type of polynomial relation appears in such low degree (compared
to the genus).


\bigskip

\setcounter{equation}{0}
\section{Geometric relations}
\bigskip

In this section we show how the simple geometric idea outlined in the
introduction produces
several families of nontrivial relations in the tautological ring
$R^*(\M_g)$. In fact,  as explained in the Appendix, 
these relations  are coming from relations in $R^*(\cc_g)$, which in
turn are restrictions of relations in the Chow ring $A^*(\ov\cc_g)$. 

\bigskip

For each $d\ge 2$, let $\cg$ be the $d$-fold fiber product of ${\cal
C}_g$ over $\M_g$, and let $p:\cg\ra \M_g$ be the forgetful
map. Denote by $D_{ij}$ the diagonal $x_i=x_j$, and let $\psi_j\in
A^1(\cg)$ be the {\em pullback} of $\psi$ from $\cc_g$ along the
map that forgets all the marked points except $x_j$. (We found it more
convenient for this paper to use this pullback class rather then the
customary  first chern class of the relative cotangent 
bundle at the marked point $x_j$.) 
Consider  
\bear
y_{d,g}=\{\; (C,x_1, \dots, x_d)\in \cg\;|\;
h^0(\O_C(x_1+\ldots+ x_d))\ge 2\;\}, 
\label{def.y.cl}
\eear
the locus corresponding to degree $d$, genus $g$ covers of a
non-rigid $\P^1$ which have all the points in a fiber marked. 

\begin{theorem}
\label{T.rel.y} 
For each $d\ge 2$, $g\ge 2$ and $b\ge 0$ we have the 
following relation in the tautological ring of $\M_g$ in dimension 
$a=g+b+1-2d$:
\bear
(2d+\kappa_0) p_* (y_{d,g}\psi_1^b)=(d-1)\kappa_{b-1} 
\cdot p_* (D_{12}\cdot  y_{d,g}) 
\label{main.rel}
\eear
Furthermore, the class $y_{d,g}$ is given 
\bear
y_{d,g}=\l[\;{c_t(E^*)\over c_t(F_d)}\;\r]_{t^{g+1-d}}
\label{rel.w}
\eear
where $c_t(E^*)=c_t(E)^{-1}$ is the total chern class of 
the dual of the Hodge bundle $E$, and  
\bear
c_t(F_d)=(1-t\psi_1)\cdot (1-t\psi_2+tD_{2,1})\cdot \ldots \cdot 
(1-t\psi_d+tD_{d,1}+\dots+tD_{d,d-1})
\label{c.F.d}
\eear 
\end{theorem}
\pf Denote by  $W_{d,g}$ the moduli space of genus $g$, degree $d$
nonsingular covers of a {\em nonrigid} $\P^1$ (in other words, we 
also mod out by automorphisms of the target, not only of the domain). 
Its pushforward to $\M_g$ defines a class $w_{d,g}$ in the Chow ring 
$A^{g+2-2d}(\M_g)$.  For each  $b\ge 0$, $\kappa_{b-1} w_{d,g}$ can be 
computed by first adding a marked point to the domain and then marking
all the other $d-1$ points in the same fiber. Therefore 
\best
{1\over (d-1)!}\cdot p_*\l( y_{d,g} \psi_1^b \r)= \kappa_{b-1} w_{d,g}
\eest 
On the other hand, $w_{d,g}$ can be also computed by first marking one
of the $r=2d+\kappa_0$ ramification points of the cover and then marking the
remaining $d-2$ points in the same fiber. 
\best
{1\over (d-2)!} p_*\l(y_{d,g}\cdot D_{12}\r) =(2d+\kappa_0) w_{d,g}
\eest 
Combining the previous two displayed equations immediately gives
(\ref{main.rel}). 

The expression (\ref{rel.w}) follows by Porteous formula as explained
in  \cite{fa} (see also  Section 7 of \cite{mu}). Let 
$F_d=p_* (\O_{\cc_g^{d+1}}(\Delta_{d+1})/
\O_{\cc_g^{d+1}})$ be the jet
bundle at $d$ points, and $E^*=R^1p_*\O_{\cc_g^{d+1}}$ be
the (pullback) of the dual of the Hodge bundle. Now look at the sequence
\best
0\longra \O_{\cc_g^{d+1}} \longra 
\O_{\cc_g^{d+1}} (\Delta_{d+1}) 
\longra
\O_{\cc_g^{d+1}}(\Delta_{d+1})/ \O_{\cc_g^{d+1}} 
 \longra 0
\eest
which gives us the sequence 
\best
0\longra F_d \ma\longra^\al  E^* 
\longra
R^1p_*( \O_{\cc_g^{d+1}}(\Delta_{d+1}))
 \longra 0
\eest
The first sheaf is locally free of rank $d$, while the second is
locally free of rank $g$. 
Since 
\best
y_{d,g}=\{\; (C,D)\;|\;h^0(C,D)\ge 2\; \}=
\{\; (C,D)\;|\;  h^1(C, D)\ge g-d+1\;  \}
\eest
then $y_{d,g}$ is exactly the locus where the rank of $\al$ drops by one, 
so $y_{d,g}$ is the Chern class $c_{g+1-d}$ of the virtual
bundle $E^*-F_d$, giving (\ref{rel.w}). 

Furthermore, by the natural filtration of $F_d$, 
\best
c(F_d)=c(F_{d-1})(1-\psi_d+D_{d,1}+\dots+D_{d,d-1})
\eest
since the first chern class of the relative dualizing sheaf of 
${\cal C}^d_g\longra {\cal C}^{d-1}_g$ is 
$
\psi_d-D_{d,1}-\dots-D_{d,d-1}
$. This gives (\ref{c.F.d}). \qed

\bigskip

\begin{rem} For each $b\ge 0$, relation (\ref{main.rel}) gives a 
(homogeneous) polynomial relation among $\kappa$ classes,  whose structure we
will study in this paper.  For $b=0$, 1 the relation is respectively
\bear
&&p_* y_{d,g}=0 
\label{main.rel.1}
\\
&& (2d+\kappa_0)\cdot p_*(\psi_1\cdot y_{d,g})= (d-1)\cdot \kappa_0
\cdot p_* (D_{12}\cdot  y_{d,g}) 
\label{main.rel.2}
\eear
The relation (\ref{main.rel.2})  was studied in \cite{ziv} for low degree
$d=2,\;3$; the relations (\ref{main.rel.1}) and (\ref{main.rel.2})  
were also mentioned  in \cite{fa} as potential candidates for proving
Theorem \ref{T.Fab}. Even though it is very plausible that the
relations for $b=0$ and 1 are enough to prove Theorem \ref{T.Fab}, our
proof involves also the relations for $b\ge 2$.
\end{rem}

\medskip

Note that equation (\ref{main.rel}) from Theorem \ref{T.rel.y}
gives a polynomial relation among $\kappa$ classes in the
tautological ring $R^*(\M_g)$. But the specific form of this polynomial
looks very complicated at first glance, and it is not clear, a priori, whether this
relation is ever nontrivial. 

We will show that in fact, the relations are encoded in a simple form
by a generating function. 
\begin{defn}\label{Def.q.c} 
Consider the positive integers $q_{k,j}$ for $k\ge j\ge
0$ (and vanishing otherwise) defined recursively by the relation
\bear
q_{k,j}=(2k+4j-2)q_{k-1,j-1}+(j+1)q_{k-1,j}+\sum_{m=0}^{k-1}
\sum_{l=0}^{j-1}q_{m,l}q_{k-1-m,j-1-l}
\label{rec.eq.q}
\eear
with initial condition $q_{0,0}=1$. Next, define recursively the
coefficients $c_{k,j}$ for $k\ge 1$ and $k\ge j\ge 0$ (and vanishing
otherwise), by the relation
\bear
q_{k,j}=(2k+4j)c_{k,j}+(j+1)c_{k,j+1}
\label{q=c}
\eear
for all $k\ge 1$ and $k\ge j\ge 0$.
\end{defn}
The coefficients $c_{k,j}$ are not necessarily integers, nor are they
always positive. Notice that, for example, (\ref{q=c}) gives 
$c_{k,k}=q_{k,k}/(6k)$, and  Proposition \ref{P.gen.c.aa} gives a very simple 
generating function for $c_{k,k}$. Furthermore,  we will see that 
$c_{k,0}={ B_{k+1}\over k(k+1)}$ (where $B_k$ are the Bernoulli
numbers),  which is not at all obvious from their definition above. 
\medskip

With the notations of Definition \ref{Def.q.c}, 
\begin{theorem}
\label{T.rel.simple}
For each $g,d\ge 2$ and each $b\ge 0$, relation  (\ref{main.rel})
 gives
 the following relation in $R^{g+1+b-2d}(\M_g)$
\bear
\l[ \;\exp \l(-\sum_{a=1}^\infty  x^{a}\kappa_{a} \sum_{j=0}^a 
c_{a,j} u^j \r)\cdot
 \l(\kappa_{b-1}-2\ma\sum_{a=0}^\infty \kappa_{a+b} x^{a+1} 
\sum_{j=0}^{a} q_{a,j}u^{j+1} \r)\;\r]_{x^{g+2-2d}u^d}
=0
\label{main.rel.2.2}
\eear
Furthermore, for $b=0$, this relation simplifies to 
\bear
\l[ \;\exp\l(- \sum_{a=1}^\infty x^{a}\kappa_{a} \sum_{j=0}^{a} 
c_{a,j} u^j \r) \;\r]_{x^{g+1-2d}u^d}=0
\label{main.rel.1.2}
\eear
\end{theorem}
The proof of this theorem essentially occupies the rest of the paper.
\bigskip

An immediate consequence of Theorem \ref{T.rel.simple} is the following
\begin{theorem}[Faber's Conjecture] 
\label{T.Fab}
The $[g/3]$ classes $\kappa_1,
\dots, \kappa_{[g/3]}$ generate the tautological ring $R^*(\M_g)$, with no
relations in degrees less or equal then $[g/3]$.
\end{theorem}
\pf The fact that there are no relations in degree less or equal to
$[g/3]$ follows from Harer's Stability result \cite{ha},
as pointed out in \cite{fa}. So it is enough to show that for each
$a\ge [g/3]+1$, Theorem \ref{T.rel.simple} gives us a relation
\bear
\kappa_a=\mbox{ polynomial in }\kappa_1, \dots, \kappa_{a-1}
\label{k=lower}
\eear
For $b\ge 2$, relation (\ref{main.rel.2.2}) gives in dimension $a=g+1+b-2d$
\best
q_{a-b,d-1} \kappa_{a}+ \mbox{ polynomial in }\kappa_1, \dots, \kappa_{a-1}=0
\eest
Since $q_{a-b,d-1}$ is a positive integer as long as $a-b\ge d-1$, 
this gives a relation of type 
(\ref{k=lower}) for each $a=g+1+b-2d\ge b+d-1$. Assume $a\ge
(g+5)/3$. Then we can find a positive integer $d$ in the range 
$(g+3-a)/2\le d \le (g+2)/3$ and choose $b=a+2d-g-1\ge 2$. Therefore, 
we get a relation of type (\ref{k=lower}) for each $g\ge 2$ and $a \ge
(g+5)/3$.

To cover the range $(g+4)/3\ge a \ge [g/3]+1$, we will use the
relations for $b=0$ and 1. In dimension $a=g+1+b-2d$, relation
(\ref{main.rel.1.2}) for $b=0$ and  (\ref{main.rel.2.2}) for $b=1$
become respectively 
\best	
&c_{a,d}\kappa_a=\mbox{ polynomial in }\kappa_1, \dots, \kappa_{a-1}
\\
&(\kappa_0 c_{a,d} +2q_{a-1,d-1})\kappa_a=
\mbox{ polynomial in }\kappa_1, \dots, \kappa_{a-1}
\eest
By Lemma \ref{L.c.aa}, the coefficient of $\kappa_a$ in both these 
relations is positive for $a=d, d+1$. 
Since $a=g+1+b-2d$ (with $b=0,1$) then $a=d$ gives
$3a=g+1+b$, while $a=d+1$ gives $3a=g+3+b$. Together, they cover the cases 
$3a=g+1, \dots, g+4$, exactly all 
missing cases in the range $(g+4)/3 \ge a\ge [g/3]+1$. \qed

\bigskip
\begin{rem} We believe that the relations (\ref{main.rel.2.2}) just for 
$b=0$ and 1 would be enough to prove Theorem \ref{T.Fab}. For that,
one needs to prove that both $c_{a,d}\ne 0$ and 
$(2a-4d-6)c_{a,d} +2q_{a-1,d-1}\ne 0$  for all $a\ge d\ge 1$. Then the
previous two displayed equations give a relation of type 
(\ref{k=lower}) for each $g\ge 2$ and $a \ge [g/3]+1$. (We checked
using a computer that this is the case for $a\le 60$, but we 
don't have a proof at this moment for why this would be true in general.)

However, having relations for $b\ge 2$ makes the proof of Theorem
\ref{T.Fab} much easier. Furthermore,
notice that for each $g\ge 2 $ and $a\ge [g/3]+1$ fixed, the relations
(\ref{main.rel.2.2}) obtained for various $d$ and $b$ (with
$a=g+1-2d+b$) are linearly independent. This is because in
(\ref{main.rel.2.2}) the coefficients of $\kappa_1^{a-j} \kappa_j$
are zero for $j\le b-2$, and nonzero for $j=b-1$ (when $b\ge 1$). 
More generally, note that there are no monomials in only the
variables $\kappa_1, \dots, \kappa_{b-2}$ appearing in
(\ref{main.rel.2.2}) for $b\ge 3$.

\end{rem}

We conclude this section with a simple generating function for some of
the relations (\ref{main.rel.2.2}):
\begin{prop}
\label{P.gen.c.aa}
The coefficients $c_{k,k}$ of Definition \ref{Def.q.c} are also 
given by the equality
\best
\exp\l(\sum_{k=1}^\infty  c_{k,k} t^k \r)=\sum_{a=0}^\infty
{(6k)!\over (2k)!(3k)!} \l({t\over 72}\r)^a
\eest
Furthermore, we have the following relations in $A^a(\M_g)$:
\bear
\l[\exp\l(-\sum_{j=1}^\infty c_{j,j} \kappa_j t^j   \r)\r]_{t^a}=0 
\label{rel.F}
\eear
for $a=g/3+1$ or $(g+1)/3$, and
\bear
\l[\exp\l(-\sum_{j=1}^\infty c_{j,j} \kappa_j t^j  \r)
\l(\kappa_{b-1} t^{b-1}-2 \kappa_b t^b- 
12\sum_{j=1}^\infty j c_{j,j} \kappa_{j+b} t^{j+b} \r)
\r]_{t^a}=0 
\label{rel.top.b}
\eear
for  $a=(g-1)/3+b$ or $(g+1)/3+b$ when $b\ge 1$.
\end{prop}
The proof of this Proposition is in the last section of this paper.

\bigskip
\begin{rem} Faber mentioned in \cite{fa} that the relation (\ref{rel.F}) is
the unique (up to a scalar multiple) relation in degree $a=(g+1)/3$. He
also claimed that there is a unique relation in degree $a=(g+2)/3$; 
relation (\ref{rel.top.b}) for $b=1$ gives then its explicit form.
\end{rem}


\bigskip

\setcounter{equation}{0}

\section{Generating functions}
\bigskip

It is obvious that Theorem \ref{T.rel.y} gives
polynomial relations in $\kappa$ classes, but the specific form of
these polynomials looks very complicated at first glance. We begin by
showing that the relations are encoded in a simple form by a
universal generating function $G$.
\begin{defn}\label{Def.G} 
Let  $G(x,w)=\ma\sum_{k=0}^\infty \ma\sum_{j=0}^\infty x^k
\al_{kj}w^j$ 
be the unique formal power series in $x$ and $w$ which solves the
recursive equation 
\bear
&&xwG_{ww}=w (G_w)^2+(1-x)G_w-1
\label{G.w.ode}
\\
&&G(x,0)=-\sum_{a=2}^\infty {B_{a}\over a(a-1)} x^{a}
\label{G.ic}
\eear
where $B_a$ are the Bernoulli numbers defined by 
\bear
{x\over e^x-1}=\ma\sum_{a=0}^\infty B_a{x^a\over a!} = 1-{x\over 2}+ 
\ma\sum_{i=1}^\infty B_{2i}{x^{2i}\over (2i)!}
\label{bernoulli}
\eear 
\end{defn}

We next show that the
relations from Theorem \ref{T.rel.y} take a simple form 
in terms of this generating function $G$.
\begin{theorem}\label{T.rel} For each $g\ge 2$, $d\ge 2$ and $b\ge 0$, 
the relation (\ref{main.rel}) is equivalent to  
\bear\label{main.rel.G}
\l[ \;\exp\l({1\over t}p_* G(t\psi,w) \r) \cdot 
p_*\l( \bspa (2w G_w(t\psi,w)+1)\psi^{b}
\r)\; \r]_{t^{g+2-2d}w^{d}}=0
\eear
where $p:\cc_g\ra \M_g$ is the forgetful map.
\end{theorem}
Note that in terms of the power series expansion of $G$ 
\best
{1\over t}p_* G(t\psi,w)&=&\ma\sum_{a=0}^\infty \ma\sum_{j=0}^\infty
t^{a-1}\kappa_{a-1}  \al_{aj}w^j 
\\
p_*\l((2G_w(t\psi,w)+1)\psi^b\r)&=&\kappa_{b-1}+2\ma\sum_{a=0}^\infty
 \ma\sum_{j=1}^\infty t^{a}\kappa_{a+b}  j \al_{aj}w^j 
\eest
where $\kappa_{-1}=0$, $\kappa_0=2g-2$. Furthermore, for $b=0$ 
 relation (\ref{main.rel.G}) simplifies to 
\bear
&\l[ \;\exp\l({1\over t}p_* G(t\psi,w)\; \r) \r]_{t^{g+1-2d}w^d}=0&
\label{main.rel.1.1}
\eear

\bigskip
\subsection{Proof of Theorem  \ref{T.rel}} 
\bigskip

The proof of Theorem \ref{T.rel} is done is several steps. The first
step is to expand (\ref{rel.w}). We begin with  some notations. 
If $1\le i_1<\dots <i_k\le d$ is a sequence of
integers, we denote by
\best
D_{i_1,\dots, i_k}
\eest
the (closed) stratum of $\cg$ where all the points $x_{i_l}$ are equal for
$l=1,\dots, k$. More generally,  the strata of $\cg$ are in one to one 
correspondence with partitions of the set $\{x_1,\dots, x_d\}$. 
Given an unordered partition 
$\{ J_1,\dots J_k\}$ of $\{ x_1,\dots, x_d\}$ we denote by 
\best
\De_{J_1,\dots, J_k}
\eest
the codimension $d-k$ multi diagonal in $\cg$  where all points in each 
$J_i$ are equal.  Given such a strata $\De_{J_1,\dots, J_k}$ we will denote by
$x_{*i}$ any one of the points of $J_i$ and by $l_i=\ell(J_i)$ the
number of points in $J_i$.  

In what follows, we will denote by $X_d(t)=(c_t(F_d))^{-1}$ and 
\bear
Y_d(t)&=&{p^* c_t(E^*)\over c_t(F_d)}=p^*c_t(E^*)\cdot X_d(t)
\label{def.Y.d}
\eear
both thought as elements of $A^*(\cg,\Q[t])$, where $F_d$ and $E$ are
defined in Theorem \ref{T.rel.y}.

\begin{prop} Using the notations above,
\bear
{p^*c_t(E^*)\over c_t(F_d)}= 
p^* \exp\l(\sum_{a=1}^\infty {B_{a+1}\kappa_a t^a\over a(a+1)}\r)
\cdot 
\sum_{r=0}^\infty\;
\sum_{\{J_1,\dots, J_r\}} t^{d-r}\De_{J_1,\dots, J_r} 
\prod_{i=1}^r G_{\ell(J_i)}(t\psi_{*i})
\label{rel.w.2}
\eear 
where the last sum is over all partitions $\{ J_1,\dots J_r\}$ of 
$\{ x_1,\dots, x_d\}$, and the formal power series 
\bear
G(x,w)= \ma\sum_{d=0}^\infty G_d(x) {w^d\over d!}= 
\ma\sum_{k=0}^\infty \ma \sum_{j=0}^\infty  \al_{kj} x^k w^j
\label{exp.ser.G.2}
\eear
satisfies the relation (\ref{G.w.ode}) with initial condition  (\ref{G.ic}).
\label{L.1st.expr}
\end{prop}
\pf  Mumford proved in
\cite{mu} that on $\ov\M_g$ the Chern character of the Hodge bundle is   
\bear
ch(E_g)=g+\ma\sum_{i=1}^\infty {B_{2i}\over (2i)!}\l(\kappa_{2i-1}+
{1\over 2}\xi_*\ma\sum_{m=0}^{2i-2}(-1)^m \phi_1^m\phi_2^{2i-2-m}   \r)
\label {mu.form}
\eear
But then the general formula expressing the total chern class in terms of 
the components of the chern character
\best
\sum_{n=0}^\infty c_n t^n=
\exp \l( \sum_{j=0}^\infty (-1)^{j-1} (j-1)! \cdot ch_j t^j\r)
\eest
(see for example \cite{fa}) implies that on $\M_g$
\bear
c_t(E_g^*)=
\exp \l(- \sum_{i=1}^\infty  {B_{2i}\cdot \kappa_{2i-1}t^{2i-1}
\over 2i(2i-1)}\r) 
\label{form.c.t.E}
\eear
This accounts for the first factor in (\ref{rel.w.2}). 

Next, for each $i=1,\dots, d$ we focus on the factor 
$(1-t\psi_i+tD_{i,1}+\dots+tD_{i,i-1})^{-1}$. Expanding in power series,
this is equal to   
\best
\sum_{k=0}^\infty t^k
(\psi_i-D_{i,1}-\dots-D_{i,i-1})^k
=\sum_{k=0}^\infty (-t)^k {k \choose a,a_1,\dots, a_{i-1}}
(-\psi_i)^{a}(D_{i,1})^{a_1}\dots (D_{i,i-1})^{a_{i-1}}
\eest
But on $\cg$ $D_{ij}^2=-\psi_i D_{ij}$ so if $a>0$ 
$D_{ij}^a=(-\psi_i)^{a-1} D_{ij}$ while 
if $1\le n_1<\dots< n_l<i$ then 
$D_{i,n_1}\cdot\ldots \cdot D_{i,n_l}=D_{i,n_1,\dots,n_l}$. Thus 
separating $a_n>0$ from $a_n=0$ for $n=1,\dots,i-1$ in the previous
displayed equation we get
\best
\sum_{k=0}^\infty (-t)^k \sum 
{k \choose a,a_{n_1},\dots, a_{n_l}}
(-\psi_i)^{a+\ma\sum_{m=1}^l (a_{n_m}-1)}D_{i,n_1,n_2,\dots, n_l}
\eest
where the second sum is over all subsets $N=\{n_1,\dots n_l\} \subset 
\{ 1,\dots i-1\}$ (including the empty subset) and all $a\ge 0$, 
$a_{n_1},\dots a_{n_l}\ge 1$. Let $F_l$ be the power series in $x$:
\best
F_l(x)=\sum_{a\ge 0}\sum_{a_m\ge 1} 
{k \choose a,a_1,\dots, a_l}x^{k-l} (-1)^l
\eest 
Then 
\bear
\frac1{1-t\psi_i+tD_{i,1}+\dots+tD_{i,i-1}}&=&
\sum_N \; F_{\ell(N)}(t\psi_i) \cdot  D_{i,N} t^{i-\ell(N)-1}
\label{psi.i}
\eear
where the sum is over all subsets $N$ of $\{ 1, \dots, i-1\}$ and 
$\ell(N)$ is the number of elements of $N$.

Finally, by (\ref{c.F.d}),
\bear
X_d(t)=
\frac1{1-t\psi_1}\cdot \frac1{1-t\psi_2+tD_{2,1}}\cdot \ldots \cdot 
\frac1{1-t\psi_d+tD_{d,1}+\dots+tD_{d,d-1}}
\label{def.xd}
\eear
But relation (\ref{psi.i}), applied for each $i=1,\dots, d$, then gives 
\best
X_d(t)=\sum D_1D_{2,N_2}\dots D_{d,N_d} 
\prod_{i=1}^d F_{\ell(N_i)}(t\psi_i)
\eest
where the sum is over all subsets $N_k$ of $\{1,\dots,k-1\}$ for each
$k=1,\dots, d$. Now, by definition, $F_d$ is symmetric in $x_1,\dots,
x_d$, which combines with the previous relation to give an expression 
\bear 
X_d(t)=\sum_{r=1}^d\sum_{J_1,\dots J_r} 
\Delta_{J_1,\dots, J_r}t^{d-r} \prod_{j=1}^ r  G_{l_j}(t\psi_{*j})
\label{prod.psi}
\eear
where $J_1,\dots, J_r$ is now a partition of $\{1,\dots, d\}$; 
the function $G_l(x)$, to be determined, is a power series in $x$
which depends only on the length $l\ge 1$. Note that the codimension of 
$\Delta_{J_1,\dots, J_r}$ is $d-r$, which accounts for the power of
$t$. 

Combining relations (\ref{form.c.t.E}) and (\ref{prod.psi}) gives (\ref{rel.w.2}). 
\bigskip

The final step in the proof of Proposition \ref{L.1st.expr} 
is to prove the recursive formula (\ref{G.w.ode}) 
for the generating function 
$G=G_0(x)+\ma\sum_{l=1}^\infty G_l(x)  w^l/l!$, where $G_0$ is chosen so
that the initial condition (\ref{G.ic}) is satisfied. Let 
$\pi:{\cal C}^{d}_g\ra {\cal C}^{d-1}_g$ be the map that forgets
$x_d$. Then by the
definition (\ref{def.xd}) of $X_d$, for any $d\ge 1$ we have the
relation 
\bear
\pi^*(X_{d-1}(t))= X_d(t)\cdot (1-t\psi_d+tD_{d,1}+\dots + tD_{d,d-1})
\label{rec.rel.xd}
\eear

Note that if $J'$ is a partition of $\{1,\dots, d-1\}$ then the
pullback $\pi^*(\De_{J'})=\De_{J',\{d\}}$. Moreover, if $J_1,\dots,
J_r$ is a partition of $\{1, \dots, d\}$ such
that $d\in J_r$ then 
\[
\Delta_{J_1,\dots, J_r}\cdot D_{d,i}=\l\{
\begin{array}{ll} 
-\psi_{*r}\; \Delta_{J_1,\dots, J_r} &\mbox{ if } i\in J_r
\\
\Delta_{J_1,\dots,\wh J_m, \dots, J_m\cup J_r}&\mbox{ if } i\in J_m,
\mbox{ for } m\ne r 
\end{array}
\r.
\]
so
\best
\Delta_{J_1,\dots, J_r}\l(1-\psi_d+\sum_{j=1}^{d-1} D_{d,j}\r)
=\Delta_{J_1,\dots, J_r}-l_r\Delta_{J_1,\dots, J_r}\psi_{*r}+
\sum_{m=1}^{r-1}l_m \Delta_{J_1,\dots,\wh J_m, \dots, J_m\cup J_r} 
\eest
We plug this and relation (\ref{prod.psi}) into (\ref{rec.rel.xd}) and 
collect the coefficient of the stratum $\Delta_{J_0}$ on both sides, where 
$J_0=\{1,\dots,d\}$. We get
\bear
\delta_{d,1}=G_{d}(x)-dxG_{d}(x)+\sum_{l_1=1}^{d-1} {d-1 \choose l_1}  l_1
G_{l_1}(x)G_{d-l_1}(x)
\eear  
for all $d\ge 1$, where $x=\psi_{*1}=\psi$ and $\delta_{d,1}$ is the
Kronecker symbol. Multiplying with $w^{d-1}/(d-1)!$ and summing
over $d$ gives 
\best
1=G_w-x(wG_{w})_w+w(G_w)^2
\eest
which is equivalent to  (\ref{G.w.ode}), thus completing 
the proof of Proposition \ref{L.1st.expr}.\qed

\bigskip
The final step in the proof of Theorem \ref{T.rel} is to use 
Proposition \ref{L.1st.expr} to express the relations
(\ref{main.rel}) in terms of the universal
generating function $G$. 

\begin{lemma} \label{L.rel.2}
The relation (\ref{main.rel}) is equivalent to the following
\bear
&\l[ \;\exp\l({1 \over t} p_* G(t\psi,w)\r)\cdot 
p_*\l( \bspa (2w G_w(t\psi,w)+1)\psi^b\r) \;\r]_{t^{g+2-2d}w^d}=0&
\label{rel.rep}
\eear
\end{lemma}
\pf With the notations in this section, relation
(\ref{main.rel}) becomes
\bear
(2d+\kappa_0)\cdot [\;p_*(Y_d \cdot \psi_d^b)\;]_{t^{g+1-d}}
=(d-1)\kappa_{b-1}\cdot [\;p_*(Y_d \cdot D_{d,1}) \;]_{t^{g+1-d}}
\label{rel.3}
\eear
We first show that after some algebraic manipulation,
(\ref{rel.3}) is equivalent to
\bear
2d \cdot  [\;p_*(Y_d \cdot \psi_d^b)\;]_{t^{g+1-d}}
=-\kappa_{b-1}\cdot [\;p_*(Y_d)\; ]_{t^{g+2-d}}
\label{rel.3.3}
\eear
On one hand, relation (\ref{rec.rel.xd}) implies that for
$d\ge 2$ 
\best
0=p_*\pi^*Y_{d-1}=p_*Y_d-t p_*(\psi_d\cdot Y_d)+ (d-1)tp_*(D_{d,1}\cdot Y_d)
\eest
In particular,
\best
(d-1)[\;p_*(D_{d,1}\cdot Y_d)\; ]_{t^{g+1-d}}=
[\;p_*(\psi_d\cdot Y_d)\; ]_{t^{g+1-d}}
 -[\;p_*Y_d \; ]_{t^{g+2-d}}
\eest
On the other hand, relation (\ref{rel.3}) for $b=1$ gives
\best
\kappa_0(d-1)[\;p_*(D_{d,1}\cdot Y_d)\; ]_{t^{g+1-d}}
=(\kappa_0+2d)[\;p_*(\psi_d\cdot Y_d)\; ]_{t^{g+1-d}}
\eest
We can then solve for $[\; p_*(D_{d,1}\cdot Y_d)\; ]_{t^{g+1-d}}$ out of the 
previous two relations
\best
2d(d-1)[\; p_*(D_{d,1}\cdot Y_d)\; ]_{t^{g+1-d}}
 =-(\kappa_0+2d) [\; p_*Y_d\; ]_{t^{g+2-d}}
\eest
Plugging this back in (\ref{rel.3}) gives (\ref{rel.3.3}). 
The push-forwards of  the terms of (\ref{rel.3.3}) are given by 
Lemma \ref{push.Y} below:
\best
\kappa_{b-1} [\;p_*(Y_d)\;]_{t^{g+2-d}}&=&\kappa_{b-1} d! \l[ \exp
\l({1\over t} \;p_*G(t\psi,w) \r)  \r]_{t^{g+2-2d}w^d}
\\
2d \cdot  [\;p_*(Y_d \cdot \psi_d^b)\;]_{t^{g+1-d}}&=& 2 d! \l[ \exp
\l({1\over t} \;p_*G(t\psi,w) \r) p_*(G_w(t\psi,w)\psi^b)
 \r]_{t^{g+2-2d}w^{d-1}}
\eest
Substituting in (\ref{rel.3.3}) and dividing by $d!$ gives (\ref{rel.rep}). 
\qed

\begin{lemma}\label{push.Y} In terms of the generating function $G$, 
 \bear
\sum_{d=1}^\infty {w^dt^{-d}\over d!} \;p_* Y_d(t) &=& \exp
\l({1\over t} \;p_*G(t\psi,w) \r) 
\label{expr.1.Y.d}
\\
\sum_{d=1}^\infty {w^{d-1}t^{-d}\over (d-1)!} \;p_* (Y_d(t) \psi^b) &=& \exp
\l({1\over t} \; p_*G(t\psi,w) \r)  \cdot {1\over t} p_*( G_w(t\psi,w) \psi^b) 
\label{expr.2.Y.d}
\eear
\end{lemma}
\pf  The expressions for both $Y_d(t)$ and $X_d(t)$ are given by Proposition
\ref{L.1st.expr}. For $X_d(t)$, the right hand side of
(\ref{prod.psi}), after being pushed forward by $p$, depends only on the
lengths $l_j$ of sets $J_j$. So 
\best
t^{-d}p_*X_d(t)=\sum_{r}\sum_{l_1,\dots, l_r\ge 1} 
t^{-r} {d\choose l_1,\dots, l_r} 
\prod_{j=1}^r p_*G_{l_j}(t\psi)= d! \sum_{m} 
\prod_{l} \l( {p_*G_{l}(t\psi)\over l! t}\r)^{m_l}
\eest
where the last sum is over all partitions $m=(m_1,m_2,\dots)$ of $d$
into $m_l$ sets of order $l\ge 1$, so $\ma\sum_{l=1}^\infty  lm_l
=d$.  Multiplying both sides with $w^d$ and summing over
$d\ge 1$ we get 
\best
\sum_{d=1}^\infty {w^dt^{-d}\over d!} \;p_* X_d(t)= 
\exp \l(\sum_{l=1}^\infty 
{ w^l\over l!}\;{p_*G_l (t\psi)\over t}\r)  
\eest
But $Y_d(t)=p^*c_t(E^*) \cdot X_d(t)$ and
(\ref{form.c.t.E}) combined with (\ref{G.ic}) implies that 
\best
c_t(E^*)=\exp\l({p_*G_0(t\psi)\over t} \r)
\eest
Therefore
\best
\sum_{d=1}^\infty {w^dt^{-d}\over d!} \;p_* Y_d(t)=  \exp \l( {1\over t}
\sum_{l=0}^\infty  {w^l\over l!} \;p_*G_l (t\psi)\r)  
\eest
which gives (\ref{expr.1.Y.d}).

Similarly, if we assume that $d\in J_1$ then 
\best
t^{-d}p_*(X_d(t)\cdot \psi_d^b)&=&\sum_{r}\sum_{l_1,\dots, l_r\ge 1} 
t^{-r} {d-1\choose l_1-1,\dots, l_r} p_*(G_{l_1}(t\psi_{*1})\cdot
\psi_{*1}^b) 
\prod_{j=2}^r p_*G_{l_j}(t\psi)
\\
&=& (d-1)! 
\sum_{r}\sum_{l_1,\dots, l_r\ge 1} \l(  {p_*(G_{l_1}(t\psi_{*1})\cdot
\psi_{*1}^b) \over (l_1-1)! t}\r) \prod_{j=2}^r {p_*G_{l_j}(t\psi)\over
l_j! t}
\eest
Therefore after multiplying by $c_t(E^*)w^{d-1}/(d-1)!$ and summing over $d$
we similarly get
\bear
\ma\sum_{d=1}^\infty { w^{d-1}t^{-d}\over (d-1)!} \;p_*(\psi^b\cdot Y_d) =
\exp\l({p_* G(t\psi,w)\over t}\r){p_*( G_w(t\psi,w)\psi^b)\over t}
\label{expr.Y.3}
\eear 
completing the proof of (\ref{expr.2.Y.d}). \qed


\bigskip

\setcounter{equation}{0}
\section{Properties of the generating function $G$}
\bigskip

In order to prove Theorem \ref{T.rel.simple}, we need to better understand the
structure of the generating function $G$ appearing in Theorem
\ref{T.rel}. We begin with the following
\begin{lemma} Equation (\ref{G.w.ode}) determines the function
$G_w$ recursively as a power series in  $x$:
\bear
G_w(x,w)&=&{-1+\sqrt{1+4w}\over 2w}+ {x\over 1+4w}+\ma\sum_{k=1}^\infty 
\ma\sum_{j=0}^k x^{k+1} q_{k,j}(-w)^j(1+4w)^{-j-{k\over 2}-1}
\label{expans.G.w.2}  
\eear
where the coefficients $q_{k,j}$ are those appearing in Definition \ref{Def.q.c}.
\label{L.sol.G.w}
\end{lemma}
\pf Denote   
\bear
G_w(x,w)=\ma\sum_{k=0}^\infty x^k A_k(w)
\label{G.w=A}
\eear
where $A_k(w)$ is a function of $w$. If we let $x=0$ in the 
equation (\ref{G.w.ode}), we get a quadratic equation
\best
wA_0(w)^2+A_0(w)-1=0
\eest
with solution $A_0(w)={-1\pm \sqrt{1+4w}\over 2w}$. Since $A_0(w)$ is a
power series in $w$ then
\bear
G_w(0,w)=A_0(w)={-1+\sqrt{1+4w}\over 2w}
\label{form.A.0}
\eear 
which accounts for the first term in (\ref{expans.G.w.2}).

We next separate $A_0(w)$ from $G_w$. For that we write 
\bear
G_w(x,w)=A_0(w)+xZ(x,w) \qquad \mbox{ where }\quad  
Z(x,w)=\sum_{k=0}^\infty x^{k} A_{k+1}(w)
\label{G.w=Z}
\eear 
and plug in equation (\ref{G.w.ode}), using the fact that $A_0(w)$ is 
a solution for $x=0$ (so there will be no free term in $x$):
\best
x^2wZ_w+xwA_0'=x^2w Z^2+2xwA_0Z+(1-x)xZ-xA_0
\eest
Dividing by $x$ and rearranging we get
\best
(wA_0)'+xwZ_w -xw Z^2+xZ=Z(1+2wA_0)
\eest
But $1+2wA_0=(1+4w)^{1/2}$ while $(wA_0)'=(1+4w)^{-1/2}$ so after
multiplying with $(1+4w)^{-1/2}$ the equation becomes
\bear
Z=(1+4w)^{-1}+x(1+4w)^{-1/2}(wZ_w-wZ^2+Z)
\label{eq.Z}
\eear
In particular we see that $Z(0,w)=(1+4w)^{-1}$ which combined with 
(\ref{G.w=Z}) explains the second term in (\ref{expans.G.w.2}). 

To prove the remaining part of the expansion (\ref{expans.G.w.2}), we
let $Q=(1+4w)Z$ and make a change of variable  
\bear
u={-w\over 1+4w}\quad \mbox{ so } \quad (1+4w)^{-1}=1+4u \quad \mbox{
 and } \quad  w={-u\over 1+4u}
\label{ch.uw}
\eear
Then $Z=(1+4u)Q$ while $(wZ)_w=-(1+4u)^2\cdot (wZ)_u=(1+4u)^2
(uQ)_u$. Plugging into (\ref{eq.Z}) and  dividing both sides by $1+4u$
we get 
\bear
Q=1+x(1+4u)^{1/2}\l((1+4u)(uQ)_u+uQ^2\r)
\label{eq.Q.1}
\eear
It is clear then by induction on powers of $x$ that $Q$ has a formal power
series expansion 
\bear
Q=\sum_{k=0}^\infty\sum_{j=0}^k x^k(1+4u)^{k/2} q_{k,j}u^j 
\label{exp.Q}
\eear
with $q_{k,j}$ positive integers (and $q_{0,0}=1$). 
Since in terms of the  original variable $w$ the function $Q$ is given by 
\bear
G_w(x,w)=G_w(0,w)+x(1+4w)^{-1}Q(x,w)
\label{G.w=Q}
\eear
then the expansion (\ref{exp.Q}) implies the last part of the expansion 
(\ref{expans.G.w.2}).

It only remains to show that the coefficients $q_{k,j}$ satisfy the
relation (\ref{rec.eq.q}) from Definition \ref{Def.q.c}. But
\best
(1+4u)(uQ)_u&=&\sum_{k=0}^\infty\sum_{j=0}^k x^k(1+4u)^{k/2} 2k
q_{k,j}u^{j+1}+\sum_{k=0}^\infty\sum_{j=0}^k x^k(1+4u)^{k/2+1}
q_{k,j}(j+1)u^{j}
\\
&=&\sum_{k=0}^\infty\sum_{j=0}^k x^k(1+4u)^{k/2}u^j 
(q_{k,j-1}(2k+4j)+q_{k,j}(j+1))
\eest
after expanding $(1+4u)$ in the last term of the first displayed
equation and then re-indexing the sum. Plugging this in (\ref{eq.Q.1}) 
and collecting the coefficient of $x^k(1+4u)^{k/2}u^j$ gives 
 the recursive relation (\ref{rec.eq.q}).\qed 

\bigskip

Next we focus on the expansion for $G$

\begin{lemma} \label{L.exp.G} 
The solution $G(x,w)$ of the equations (\ref{G.w.ode})
and (\ref{G.ic}) has the form:
\bear
G(x,w)&=&G(0,w)+{x\over 4}\ln(1+4w)-\ma\sum_{k=1}^\infty 
\ma\sum_{j=0}^\infty x^{k+1} c_{k,j}(-w)^j(1+4w)^{-j-{k\over 2}} 
\label{form.G.2}
\eear
where the coefficients $c_{k,j}$ are related to 
the coefficients $q_{k,j}$  via the formula (\ref{q=c}).  
\end{lemma} 
\pf After the change of variable (\ref{ch.uw}) relation 
(\ref{G.w=Q}) combined with the fact that $G_w=-(1+4u)^2 G_u$ implies
that 
\bear
-G_u(x,u)=-G_u(0,u)+ x(1+4u)^{-1}Q(x,u)
\label{eq.G.u}
\eear
where $Q$ has the expansion (\ref{exp.Q}). Therefore integrating 
\best
(1+4u)^{-1} Q(x,u)={1\over 1+4u}+\ma\sum_{k=1}^\infty x^k
(1+4u)^{k/2-1} \ma\sum_{j=0}^k  q_{kj} u^j
\eest 
we get an expansion 
\bear\nonumber
-G(x,u)&=&-G(0,u)+{x\over 4}\ln(1+4u)+\ma\sum_{k=1}^\infty x^{k+1}
\l( \ma\sum_{j=0}^k(1+4u)^{k/2}   \wt c_{k,j}u^j + \al_k\r)\hskip.1in
\\
&=&-G(0,u)+{x\over 4}\ln(1+4u)+\ma\sum_{k=1}^\infty x^{k+1}
(1+4u)^{k/2}\ma\sum_{j=0}^\infty c_{k,j}u^j
\hskip.1in
\label{G+al}
\eear
The coefficients $\al_k$ are determined from the initial condition
(\ref{G.ic}). In fact, this initial condition implies that all
$\al_k=0$ (i.e. $c_{k,j}=0$ for all $j>k$). 
We were not able to show this directly from the equations,
but Lemma \ref{L.rel.uy} provides an indirect argument based on Harer's 
stability result \cite{ha}. In any event, the expansion (\ref{G+al}), 
after changing back to the $w$ variable becomes (\ref{form.G.2}).

Furthermore, differentiating relation (\ref{G+al}) with respect to $u$ gives
\best
-G_u(x,u)&=&{x\over 1+4u}+\ma\sum_{k=1}^\infty x^{k+1}
\ma\sum_{j=0}^\infty (1+4u)^{k/2-1} 2kc_{k,j}u^j+
\ma\sum_{k=1}^\infty x^{k+1}
\ma\sum_{j=1}^\infty (1+4u)^{k/2} jc_{k,j}u^{j-1}
\\
&=&x(1+4u)^{-1}\ma\sum_{k=1}^\infty x^{k}
\ma\sum_{j=0}^\infty (1+4u)^{k/2} u^j(2kc_{k,j}+4j
c_{k,j}+(j+1)c_{k,j+1}) 
\eest
The last equality follows by expanding the factor $1+4u$ in the second 
sum of the previously displayed equation and then rearranging the
terms. Comparison with (\ref{eq.G.u}) immediately implies
(\ref{q=c}). \qed 


\bigskip

As we have seen in the proof of the previous lemmas, it is 
more convenient to make the following change of variables
\bear
u={-w\over 1+4w}\quad \mbox{ and }\quad y={x\over \sqrt{1+4w}}=x(1+4u)^{1/2}
\label{xw=yu}
\eear
\begin{lemma}\label{L.ch.var} Let $P(x,w)$ be a formal power series in $x$ and
$w$. Denote by $\wh P(y,u)$ the formal power series in $y,u$ obtained from
$P(x,w)$ after the change of variables (\ref{xw=yu}). Then
\bear
\l[P(x,w) \r]_{x^{a} w^d}= (-1)^d 
\l[(1+4u)^{{a+2d-2\over 2}} \wh P(y,u) \r]_{y^{a} u^d}
\label{ch.var}
\eear
\end{lemma}
\pf To begin with, 
\best
&&\l[P(x,w) \r]_{x^{a} }= (1+4w)^{-a/2} \l[ P(y,w) \r]_{y^{a} }
\eest
But if we make the change of variable (\ref{xw=yu}), 
$du= -(1+4w)^{-2}dw$ then:
\best
\l[ (1+4w)^{-a/2} f(w) \r]_{w^d}&=& 
\oint   {(1+4w)^{-a/2} f(w)\over w^{d+1}}dw=
(-1)^d \oint  {(1+4u)^{(a+2d-2)/2} \wh f(u)\over u^{d+1}}du
\eest 
which gives (\ref{ch.var}). \qed

\begin{lemma} In terms of the expansions of Lemmas \ref{L.sol.G.w} and 
\ref{L.exp.G}, 
 after the change 
of variables (\ref{xw=yu}), the relations  (\ref{main.rel.1.1}) 
and (\ref{main.rel.G}) become respectively
\bear
&&\l[ \;\exp\l(- \sum_{a=1}^\infty x^{a}\kappa_{a} \sum_{j=0}^\infty
c_{a,j} u^j \r) \;\r]_{x^{g+1-2d}u^d}=0
\label{main.rel.1.5}
\\
&&\l[ \;\exp \l(-\sum_{a=1}^\infty  x^{a}\kappa_{a} \sum_{j=0}^\infty 
c_{a,j} u^j \r)\cdot \l(\kappa_{b-1}-2\ma\sum_{a=0}^\infty \kappa_{a+b} x^{a+1} 
\sum_{j=0}^{a} q_{a,j}u^{j+1} \r)\;\r]_{x^{g+2-2d}u^d}
\hskip-.2in=0\hskip.2in 
\label{main.rel.2.5}
\eear
Furthermore, $c_{k,j}=0$ for all $j>k$.
\label{L.rel.uy}
\end{lemma}
\pf With the expansion (\ref{form.G.2}), 
\best
{1\over t}p_*G(t\psi,w)={\kappa_0\over 4}\ln(1+4w)-
\sum_{a=1}^\infty t^a \kappa_a \sum_{j=0}^\infty 
c_{aj}(-w)^j (1+4w)^{-j-a/2}  
\eest
so after changing the variables
\bear
\exp\l(\; {1\over t}p_*G(t\psi,w)\;\r) = (1+4u)^{-\kappa_0/4}
\exp\l(\; \sum_{a=1}^\infty y^a \kappa_a  \sum_{j=0}^\infty  c_{aj}u^j \;\r)
\label{1.term.yu}
\eear
Similarly, with  the expansion (\ref{expans.G.w.2}), 
 $p_*((2wG_w(t\psi,w)+1)\psi^b)$ becomes
\best
&&\kappa_{b-1}\sqrt{1+4w}-
2\sum_{a=0}^\infty 
t^{a+1} \kappa_{a+b} \sum_{j=0}^a q_{aj}(-w)^{j+1} (1+4w)^{-j-{a\over 2}-1}  
\\
&=&(1+4w)^{1/2}\l(\kappa_{b-1}-
2\sum_{a=0}^\infty t^{a+1} \kappa_{a+b} \sum_{j=0}^a q_{aj}(-w)^{j+1}
(1+4w)^{-j-{a+1\over 2}-1}   \r) 
\eest
so after changing the variables 
\bear
p_*((2wG_w(t\psi,w)+1)\psi^b)=
(1+4u)^{-1/2}\l(\kappa_{b-1}-
2\sum_{a=0}^\infty y^{a+1} \kappa_{a+b}\sum_{j=0}^a q_{aj} u^{j+1} \r) 
\label{2.term.yu}
\eear
Next we use the change of variables (\ref{ch.var}) in both
(\ref{main.rel.1.1}) and (\ref{main.rel.G}). By Lemma \ref{L.ch.var}  
we pick up a factor of $(1+4u)^{a+2d-2\over 2}$, where $a=g+1-2d$ 
for the first relation, and $a=g+2-2d$ for the second one. But then 
in the first case 
\best
{a+2d-2\over 2}={g-1\over 2}={\kappa_0\over 4}
\eest
while for the second $(a+2d-2)/2=\kappa_0/4+1/2$. This means that 
the powers of $(1+4u)$
cancel out in both cases, giving precisely
(\ref{main.rel.1.5}) and  (\ref{main.rel.2.5}) respectively.
\medskip

Finally,  we use Harer's stability result \cite{ha} to show that 
$c_{k,j}=0$ whenever $j>k$. With the notations in the proof of 
Lemma \ref{L.exp.G}, it is enough to show that $\al_k=0$ for all $k$. Assume by
contradiction that  $\al_a\ne 0$ for some $a\ge 1$. Then for any
$g>3a$ relation (\ref{main.rel.1.5}) provides a relation between
$\kappa$ classes in degree $a$ in which the coefficient of $\kappa_a$
is a nonzero multiple of $\al_a$, thus nonzero. But this contradicts 
Harer stability theorem, which implies that for each $a$ the ring 
\best
\ma\oplus _{i=0}^{2a} H^i(\M_g,\Q)
\eest
stabilizes for $g\ge 3a$ to  $\Q[\kappa_1,\dots,\kappa_a]$ 
with no relations. \qed

\bigskip

\non {\bf Proof of Theorem \ref{T.rel.simple}}. Together, the previous
Lemmas completely prove Theorem \ref{T.rel.simple}.  More precisely,
since $c_{k,j}=0$ for $j>k$ by Lemma \ref{L.rel.uy}, then the relation
(\ref{main.rel.2.5}) becomes exactly
(\ref{main.rel.2.2}). Furthermore, Lemma \ref{expans.G.w.2} and Lemma
\ref{L.exp.G} show that the coefficients in the relation
(\ref{main.rel.2.2}) are exactly those appearing in Definition
\ref{Def.q.c}.

\begin{rem} Note that the relation (\ref{q=c}) combined with the initial
condition for $G$ implies the following combinatorial relation for the
Bernoulli numbers
\best
\sum_{k=1}^\infty {B_{k+1}\over k(k+1)} t^{k} =
\sum_{k=1}^\infty {c_{k0}}t^k={1\over 4} 
\sum_{k=1}^\infty t^k\sum_{l=0}^k (-4)^{-l}
q_{kl}\cdot {l! (k/2-1)!\over (k/2+l)!}
\eest
\end{rem}

\bigskip

We conclude this section with the 
\smallskip

\non {\bf Proof of Proposition \ref{P.gen.c.aa}.} The formula for the
generating function for the coefficients $c_{k,k}$ is given by 
 Lemma \ref{L.c.aa} below. 

Next, if we make the change of variable $t=xu$ and $y=u^{-1}$ 
in relations (\ref{main.rel.1.2})
and (\ref{main.rel.2.2}) they respectively become after re-indexing the
sums
\best
&&\l[ \;\exp \l(-\ma\sum_{a=0}^\infty \ma\sum_{j=0}^a  \kappa_{a} 
c_{a, a-j} t^{a} y^{j} \r)\cdot
 \l(\kappa_{b-1}-2\ma\sum_{a=0}^\infty \ma\sum_{j=0}^a \kappa_{a+b}
q_{a,a-j}t^{a+1}y^{j} \r)\;\r]_{t^{g+2-2d}y^{3d-g-2}}
=0
\\
&&\l[ \;\exp\l(-\ma \sum_{a=0}^\infty \ma\sum_{j=0}^a 
\kappa_{a}
c_{a,a-j} t^{a} y^{j} \r) \;\r]_{t^{g+1-2d}y^{3d-g-1}}=0
\eest
We are interested in the coefficients of $y^0$ and $y^1$. Since 
by Lemma \ref{L.c.aa} $10 c_{a,a-1}= a c_{a,a}$ while 
$10 q_{a, a-1}=(a+1) q_{a,a}$ for $a\ge 1$ and $q_{0, 0}=1$ then 
the previous two displayed equations become respectively 
\bear
&&\l[ \;\exp \l(f(t)+tf'(t){y\over 10}+ O(y^2) \r)
 \l( g(t)+\l(tg'(t)+2\kappa_b\r) {y\over 10}+ O(y^2)
 \r)\;\r]_{t^{g+2-2d}y^{3d-g-2}}\hskip-.3in
=0\hskip.3in
\label{e.f.1}
\\
&&\l[ \;\exp\l( f(t)+tf'(t){y\over 10} + 
O(y^2) \r) \;\r]_{t^{g+1-2d}y^{3d-g-1}}=0
\label{e.f.2}
\eear
where
\best
f(t)&=& -\sum_{j=0}^\infty \kappa_{a} c_{a, a} t^a 
\qquad \mbox{and}
\\
g(t)&=&\kappa_{b-1}-2t \kappa_b-
2\ma\sum_{a=1}^\infty \kappa_{a+b} q_{a,a}t^{a+1}=
\kappa_{b-1}-2t \kappa_b- 12\ma\sum_{a=1}^\infty \kappa_{a+b} a c_{a,a}t^{a+1}
\eest
Then for $3d-g-1=0$ and $a=g+1-2d$, relation (\ref{e.f.2}) becomes 
\bear
\l[ \;\exp\l( f(t) \r) \;\r]_{t^a}=0
\label{e.f.f}
\eear
while for $3d-g-1=1$ and $a=g+1-2d$ it becomes
\best
\l[ \;\exp\l( f(t) \r)  t f'(t) \;\r]_{t^a}=0
\eest
which is of course the same as (\ref{e.f.f}) 
(up to a scalar factor of $a$). Together, they give (\ref{rel.F}). 

Similarly, after subtracting $\kappa_b/5$ times (\ref{e.f.2})
from (\ref{e.f.1}), we see that in the new relation (\ref{e.f.1})  the
coefficient of $x^1$ is up to a scalar multiple $t$ times the 
derivative in $t$
of the coefficient of $x^0$. Therefore, for both $g+2-3d=0$ and 
$g+2-3d=1$,  relation (\ref{e.f.1}) implies
\best
\l[ \;\exp\l( f(t) \r) g(t) \;\r]_{t^{g+2-2d}}=0
\eest
which gives (\ref{rel.top.b}).  \qed 

\begin{lemma}
\label{L.c.aa}
We have the following equality 
\bear
\exp\l( \ma\sum_{k=1}^\infty c_{kk}z^k \r)=\ma\sum_{k=1}^\infty
{(6k)!\over (2k)!(3k)!}\l({z\over 72}\r)^k 
\label{exp.c=!}
\eear
Moreover, each $k\ge 1$
\bear
0<q_{kk}=6kc_{k,k}, \quad 0<q_{k,k}=60c_{k,k-1}
\label{q=k.c}
\eear
\end{lemma}
\pf Let 
\bear
Q_0(z)=\sum_{k=1}^\infty q_{k,k-i}z^{k+1-i}\quad \mbox{ and }\quad
Q_i(z)=\sum_{k=i}^\infty q_{k,k-i}z^{k+1-i} \quad \mbox{ for } i\ge 1
\eear 
With this notation, 
\best
uQ=\ma\sum_{k\ge i\ge 0} q_{k, k-i} x^k (1+4u)^{k/2} u^{k-i+1}=u+
\ma\sum_{i=0}^\infty x^{i-1}(1+4u)^{(i-1)/2}Q_i( xu\sqrt{1+4u})
\eest
Next we make another change of variables
\best
z=xu\sqrt{1+4u}\quad \mbox{ and }\quad y=x\sqrt{1+4u}
\eest
Then using logarithmic differentiation
\best
y_u={2 y\over 1+4u} \quad \mbox{ and }\quad 
z_u={z\over u}+{2 z\over 1+4u}={z(1+6u)\over u(1+4u)} ={y+6z\over 1+4u}
\eest
In terms of this new variables,
\best
uQ&=&u+\ma\sum_{i=0}^\infty y^{i-1}Q_i(z)
\\
(1+4u)(uQ)_u&=&1+\ma\sum_{i=0}^\infty y^{i-1}
\l[2(i-1)Q_i(z)+Q_i'(z)(y+6z)\r]
\eest
Multiplying both sides of (\ref{eq.Q.1}) by $u$ gives 
$uQ=u+[z(1+4u)(uQ)_u+ y (uQ)^2]$; substituting the
previous two displayed equations we get
\best
\ma\sum_{i=0}^\infty y^{i-1}Q_i(z)=z\ma\sum_{i=0}^\infty y^{i-1}
\l[2(i-1)Q_i(z)+Q_i'(z)(y+6z)\r]+y[z/y+\ma\sum_{i=0}^\infty y^{i-1}Q_i(z)]^2
\eest
Multiplying both sides with $y$ gives
\bear
\ma\sum_{i=0}^\infty y^{i}Q_i(z)=z\ma\sum_{i=0}^\infty y^{i}
\l[2(i-1)Q_i(z)+Q_i'(z)(y+6z)\r]+[z+\ma\sum_{i=0}^\infty y^{i}Q_i(z)]^2
\label{exp.q.k-j}
\eear
Now we can collect the coefficient of $y^0$ to get
\bear
Q_0=6z^2Q_0'+Q_0^2
\label{eq.Q.0.z}
\eear
This equation can be easily integrated, using the integrating factor
\best
P(z)=\exp\l({1\over 6} \int_0^z {Q_0(\zeta)\over \zeta^2} d\zeta \r)=
\sum_{k=0}^\infty p_k z^k
\eest
Then $P'=P {Q_0\over 6z^2}$ while 
\best
P''=P{Q_0^2\over 36z^4}+
P {Q_0'\over 6z^2}-2P {Q_0\over 6z^3}=P{Q_0^2+6z^2Q_0'-12z Q_0\over
36z^4}.
\eest 
Multiplying both sides 
of (\ref{eq.Q.0.z}) by ${P\over 6z^2}$ gives 
\best
6P'=5P+36z^2P''+72zP'
\eest
This means
\best
6(k+1)p_{k+1}=p_k(5+36k(k-1)+72 k)=(6k+1)(6k+5)p_{k}
\eest
which gives inductively
\best
p_k={(6k)!\over (3k)!(2k)!}72^{-k}
\eest
But of course since $Q_0(z)=\ma\sum_{k=1}^\infty q_{kk} z^{k+1}$ then 
\best
P(z)=\exp\l({1\over 6}\sum_{k\ge 1} q_{k,k}{z^k \over k} \r)= 
\exp\l(\sum_{k\ge 1} c_{k,k}z^{k} \r)
\eest
which gives  (\ref{exp.c=!}). 

For $i\ge 1$, the coefficient of $y^i$ in (\ref{exp.q.k-j}) gives
\best
Q_i=6z^2 Q_i'+2izQ_i+ z^2Q'_{i-1}+2Q_0Q_i+\sum_{l=1}^{i-1} Q_lQ_{i-l}
\eest 
When $i=1$ this becomes
\best
Q_1=6z^2Q_1'+2zQ_1+zQ_0'+2Q_1Q_0
\eest
It follows then that $10Q_1=Q_0'$: differentiating
(\ref{eq.Q.0.z}) gives
\best
Q_0'=6z^2Q_0''+12zQ_0'+2Q_0Q_0'
\eest
So if $F=10Q_1-Q_0'$ then subtracting the previous two displayed
equations  we get
\best
F=6z^2F'+2zF+2Q_0F
\eest
which by induction on the powers of $z$ implies $F=0$. This means that 
$10 q_{k,k-1}=(k+1) q_{k,k}$ for each $k\ge 1$. 

Taking $k=j$ in (\ref{q=c}) gives $6k c_{k,k}=q_{k,k}$, while $k=j+1$
 gives $q_{k,k-1}=(6k-4)c_{k,k-1}+kc_{k,k}$. Therefore,
 when $k\ge 1$,
\best
(6k-4)c_{k,k-1}=q_{kk-1}-{1\over 6}q_{kk}=q_{kk}\l( {k+1\over
10}-{1\over 6}\r)=q_{kk}{3k-2\over 30}
\eest
which together complete the proof of (\ref{q=k.c}). \qed
\bigskip


\bigskip
\renewcommand{\theequation}{A.\arabic{equation}}

\setcounter{equation}{0}

\section{Appendix}
\bigskip

The relations of Theorem
\ref{T.rel.y} are actually coming from relations in the tautological 
ring of $\cc_g^2$, which in turn are restrictions of relations 
in the Chow ring $A^*(\ov\M_{g,2})$, as explained below. 

We will work with the moduli space ${\cal Y}_{d,g}$ of degree $d$,
genus $g$ covers of $\P^1$ with a fixed ramification pattern over some
points of the target, as defined in \cite{i}. The arguments used
in this appendix are very similar to those in \cite{i}, to which we
refer the reader for more detailed explanations. 
 
Fix 4 distinct points $p_1,\dots p_4$ on $\P^1$ and consider the
moduli space 
\best 
{\cal Y}_{d,g}(B_{1^d}(y_1) B_{1^d}(y_2) B_{2,
1^{d-2}} B_{2, 1^{d-2}} ) 
\eest 
of genus $g$, degree $d$ relatively
stable covers of this $\P^1$ which are simply ramified over $p_3$ and
$p_4$ and have a marked point $y_1$ and $y_2$ in each fiber over $p_1$
and respectively $p_2$. The push forward via the stabilization map
$st$ defines a class in Chow ring $A^{g+3-2d}(\ov\M_{g,2})$. (Notice
that this class is not the same as the relative virtual fundamental
cycle described in \cite{gv}; however, we next restrict to covers
whose domain contains a smooth genus $g$ component, where the two
definitions agree)

Now, as in Proposition 2.8 of \cite{i}, we split the target $\P^1$ in
two ways. We use the fact that on $\ov\M_{0,4}$ the divisor corresponding to
fixing the location of the marked points $p_1, \dots, p_4$ is linear
equivalent to the boundary divisor $(p_1, p_2|p_3,p_4)$ where $p_1,
p_2$ are on a bubble and $p_3,p_4$ on the other. So the class 
$st_*
{\cal Y}_{d,g}(B_{1^d}(y_1) B_{1^d}(y_2) B_{2, 1^{d-2}} B_{2, 1^{d-2}}) $ 
can be computed by splitting the target either as $(p_1,
p_2|p_3,p_4)$ or as $(p_1, p_3|p_2,p_4)$ and using the degeneration
formula (1.23) of \cite{i}. This gives a relation in the Chow ring of
$\ov\M_{g,2}$. 

Here we are only interested in the restriction of the above relation
to the Chow ring of ${\cal C}^2_g$, so we are interested only in those
covers of $\P^1 \vee \P^1$  whose domain has a smooth genus $g$ 
component, and the rest are rational components 
(the ``symbol'' of the relation, using the terminology of \cite{i}). 

By dimension count, in the symbol, the only possibility is that on
one side we have a genus $g$, degree $d$ cover, while the cover on the other
side has only rational components, each totally ramified over the node
of the target. (Of course, the ramification pattern of the two covers has to 
match over the double point of the target.) If the genus $g$ component 
had degree less then $d$ then its push forward by $st$ would give a 
smaller dimensional class in the Chow ring. If one of the rational components
was not totally ramified over the node of the target, it would then
have at least two points in common with the degree $d$, genus $g$
component, which is impossible.
\medskip

So, for the splitting of $\P^1$ with $p_1, p_2$ on the right bubble
and $p_3, p_4$ on the left bubble, in the symbol, either  
\begin{enumerate} 
\item the rational components are on the left, each of degree 1; 
the two marked points $y_1$ and $y_2$ are either on
different components, or on the same component. After pushforward, this gives
\best
st_*{\cal Y}_{d,g}(B_{1^d}(y_1; y_2) B_{2, 1^{d-2}}  B_{2,
1^{d-2}} )+st_*{\cal Y}_{d,g}(B_{1^d}(y_1=y_2) B_{2, 1^{d-2}}  B_{2,
1^{d-2}} )
\eest
where the
notation $B_{1^d}(y_1=y_2)$ means the strata where $y_1, \;y_2$ are on
a bubble attached to a point where the cover has ramification pattern
$B_{1^d}$. Forgetting the marking of the two branch points, the last 
displayed equation becomes
\bear\label{t1}
{r(r-1)} \l( \mspa
st_*{\cal Y}_{d,g}(B_{1^d}(y_1; y_2) )+ 
st_*{\cal Y}_{d,g}(B_{1^d}(y_1=y_2) ) \r)
\eear
where $r=2d+\kappa_0$ is the total number of branch points;
\item or the rational components are  on the right, in which case 
the two marked branch points must land on rational components. There are 
two cases: they either land on the same rational component, which 
then forces that component to have degree 3, or they land on two different
degree 2 components. After pushforward, this gives respectively
\best
3st_*{\cal Y}_{d,g}(B_{1^d}(y_1) B_{1^d}(y_2) B_{3, 1^{d-3}} )
+
{2\cdot 2\over 2!} st_*
{\cal Y}_{d,g}(B_{1^d}(y_1) B_{1^d}(y_2) B_{2, 2, 1^{d-4}} )
\eest
Note that this is the pullback by the projection $\pi_{1,2}$ that forgets 
both $y_1, y_2$ of
\bear\label{t2}
{3}st_*{\cal Y}_{d,g}(B_{3, 1^{d-3}} )
+
{2} st_*
{\cal Y}_{d,g}( B_{2, 2, 1^{d-4}} )
\eear
\end{enumerate}
\bigskip

On the other hand, when we  split the target $\P^1$ with $p_1, p_3$ on 
the right bubble and $p_2, p_4$ on the left bubble, we get something 
symmetric in  $y_1, \; y_2$. So it is enough to 
consider the case where the rational components are on the left; the
other case is obtained by switching  $y_1$ with $y_2$. In the first case, the
rational  component containing the branch point has degree 2, and the 
others have degree 1. The marked point $y_1$ can be either on the degree 2
component or on a degree 1 component. After pushforward, this gives 
\best
{2} st_*{\cal Y}_{d,g}(B_{2,1^{d-2}}(y_1) B_{1^d}(y_2)  B_{2,
1^{d-2}} ) +
st_*{\cal Y}_{d,g}(B_{1,2,1^{d-3}}(y_1) B_{1^d}(y_2)  B_{2,
1^{d-2}} )
\eest
Note that this is equal to the pulled back by $\pi_2$ 
(the map that forgets $y_2$) of 
\bear\label{t3}
(r-1)\l( 2 st_*{\cal Y}_{d,g}(B_{2,1^{d-2}}(y_1) ) 
+st_*{\cal Y}_{d,g}(B_{1,2,1^{d-3}}(y_1) )\r)
\eear
\medskip

Assembling (\ref{t1}), (\ref{t2}), (\ref{t3}) as well as (\ref{t3}) in 
which $y_1$ is switched with $y_2$  we get the following relation in 
$R^{g+3-2d}({\cal C}^2_g)$:
\bear\nonumber
&&r(r-1)\l( \mspa 
{st_*{\cal Y}_{d,g}(B_{1^d}(y_1; y_2) )} + 
{st_*{\cal Y}_{d,g}(B_{1^d}(y_1=y_2) )} \r) 
\\ \nonumber
&&\hskip.5in + \pi_{1,2}^* \l( \mspa 
3 st_*{\cal Y}_{d,g}(B_{3, 1^{d-3}} )
+2 st_* {\cal Y}_{d,g}( B_{2, 2, 1^{d-4}} )\r) 
\\ \nonumber
=&&(r-1)\pi_2^* \l( \mspa   2st_*{\cal Y}_{d,g}(B_{2,1^{d-2}}(y_1) ) +
 st_*{\cal Y}_{d,g}(B_{1,2,1^{d-3}}(y_1) ) \r)
\\
&&\hskip.5in +
(r-1)\pi_1^* \l( \mspa   2 st_*{\cal Y}_{d,g}(B_{2,1^{d-2}}(y_2) ) +
 st_*{\cal Y}_{d,g}(B_{1,2,1^{d-3}}(y_2) ) \r)
\label{eq.m.gen}
\eear
This relation can be rewritten in terms of the class $y_{d,g}$ of
Theorem \ref{T.rel.y}. For that, we introduce $d$ more marked 
points $x_1, \dots, x_d$, and let $p:\cc_g^{d+2}\ra \cc_g^2$
be the map that forgets these $d$ points (we use the letter $\pi$ 
for the map that forgets the points $y_1$ or $y_2$). 
Let $y_{d,g}$ denote the class defined by
(\ref{def.y.cl}), as well as its pullback to $\cc_g^{d+2}$. 

Each term of relation (\ref{eq.m.gen}) corresponds to covers of a 
nonrigid $\P^1$ with a fixed ramification pattern over a point in
$\P^1$. Such class  can be expressed in terms of
$y_{d,g}$ by marking all the points in that fiber. For example,
\best
st_*{\cal Y}_{d,g}(B_{1^d}(y_1) )={1\over (d-1)!}
 st_*{\cal Y}_{d,g}(B_{1^d}(y_1;x_2, \dots, x_d) ) ={1\over (d-1)!}
 p_*(y_{d,g}\cdot D_{y_1, x_1})
\eest
where $D_{y_1,x_2}$ denotes the diagonal $y_1=x_1$. So (\ref{eq.m.gen})
gives
\begin{prop} With the notations above, the following relation holds in 
$R^{g+3-2d}({\cal C}^2_g)$:
 \bear\nonumber
&&r(r-1)\l( \mspa {1\over (d-2)!}p_*( y_{d,g}\cdot D_{y_1, x_1}D_{y_2,x_2})
+ {1\over (d-1)!}p_*( y_{d,g}\cdot D_{x_1, y_1, y_2}) \r) 
\\ \nonumber
&&\hskip.5in + \pi_{1,2}^* \l( \mspa 
3  {1\over (d-3)!}p_*( y_{d,g}\cdot D_{x_1, x_2, x_3})
+2 {1\over (d-4)!} p_*( y_{d,g}\cdot D_{x_1, x_2}D_{x_3,x_4}) \r) 
\\ \nonumber
&=&(r-1)\pi_2^* \l( \mspa   2{1\over (d-2)!}
p_*( y_{d,g}\cdot D_{x_1,x_2,y_1})) +{1\over (d-3)!}
p_*( y_{d,g} \cdot D_{x_1, y_1}D_{x_2,x_3}) ) \r)
\\
&&\hskip.5in +
(r-1)\pi_1^* \l( \mspa   2{1\over (d-2)!}
p_*( y_{d,g}\cdot D_{x_1,x_2,y_2})) +{1\over (d-3)!}
p_*( y_{d,g} \cdot D_{x_1, y_2}D_{x_2,x_3}) ) \r)\hskip.3in 
\label{r.y.2.2}
\eear
for all $g,d \ge 2$.
\end{prop}
Each term in the equation above can be expressed in terms of the 
generating function $G$ just as in the proof of Lemma \ref{push.Y}. 
In fact, if for each fixed $j$ we denote by 
$p_{j}$ the map that forgets all the marked points except $x_1,\dots,
x_j$, 
then the generating function of $p_{j*}y_{d,g}$ is 
\bear
\sum_{d} p_{j*}y_{d,g} {w^{d-j}t^{-d}\over (d-j)!}= 
\exp\l({1\over t} p_*G(t\psi, w) \r) \sum \Delta_{J_1,
\dots, J_r} t^{-r}\prod_{i=1}^r 
\l({\partial \over \partial w}\r)^{l_i}G (t\psi_{*i}, w) 
\label{push.Y.gen}
\eear
where the last sum is over all partitions 
$\{J_1,\dots, J_r\}$ of $\{ x_1, \dots, x_j\}$. For example, when $j=1$,
this simply becomes
\best
\exp\l({1\over t} p_*G(t\psi, w) \r) t^{-1} G_w (t\psi_{1}, w) 
\eest
while for $j=2$ it becomes
\bear
\exp\l({1\over t} p_*G(t\psi, w) \r) 
\l( t^{-2} G_w (t\psi_{1}, w)  G_w (t\psi_{2}, w)+ 
t^{-1} G_{ww} (t\psi_{1}, w) D_{1,2}\r)
\label{G.2}
\eear
Note that in such expressions, the factor 
$\exp\l({1\over t} p_*(G(t\psi, w) \r) $ is pulled
back from $\M_g$. 

Thus for example the first term of (\ref{r.y.2.2}) is nothing but
$r(r-1)$ times the
coefficient of $t^{g+2-2d} w^{d-2}$ in the expression (\ref{G.2}). 
Using (\ref{push.Y.gen}), relation (\ref{r.y.2.2}) then becomes a fairly 
complicated, but explicit relation in $A^*(\cc_g^2)$  involving
up to 4 derivatives in $w$ of the generating function $G$. One could
use the ODE (\ref{G.w.ode}) satisfied by $G$ to reduce the relation
above to a 
not much simpler one
involving, besides the factor $\exp\l({1\over t} p_*G(t\psi, w) \r)$,
an explicit polynomial in $G_w(t\psi_1,w)$, $G_w(t\psi_2,w)$ and 
$D_{1,2}$.

\bigskip
\begin{rem} 
If we  pushforward relation (\ref{eq.m.gen}) by  $\pi_2$ and 
divide by $(r-1)d$  we simply get
\best
r\cdot 
st_*{\cal Y}_{d,g}(B_{1^d}(y_1) )
=\pi_1^* st_*{\cal Y}_{d,g}(B_{2,1^{d-2}} ) 
\eest
which as before is equivalent to
\best
{1\over (d-1)!}
 p_*(y_{d,g}\cdot D_{y_1, x_1})= {1\over
(d-2)!} p_* (y_{d,g}\cdot D_{x_1, x_2})
\eest
Thus in $R^{g+2-2d}(\cc_g)$ we have the relation 
\bear
r \cdot p_* ( y_{d,g}\cdot D_{y_1, x_1})=(d-1) p_*(y_{d,g}\cdot D_{x_1, x_2})
\label{rel.c}
\eear
where $y_{d,g}$ defined by (\ref{def.y.cl}). Multiplying by $\psi_1^b$ 
and the pushing forward by $\pi_1$ gives exactly the relations
(\ref{main.rel}) of Theorem \ref{T.rel.y}. It is easy to modify the
arguments in Section 2 to see that in 
terms of the generating function $G$, 
(\ref{rel.c}) becomes
\best
\l[ \exp\l({1\over t} p_*G(t\psi, w) \r) (2w G_w(t\psi_1, w)-1)
\r]_{t^{g+2-2d}w^d} =0 
\eest
which in turn becomes relation (\ref{main.rel.c}) mentioned in the 
introduction. 
\end{rem}



\begin{thebibliography}{}

\bibitem[F]{fa} C. Faber, {\em A conjectural description of the 
tautological ring of the moduli space of curves}, Moduli of
curves and abelian varieties, 109-129, Aspects Math., E33, Vieweg,
Braunschweig, 1999. 

\bibitem[GV]{gv} T. Graber, R. Vakil,  
{\em Relative virtual localization 
and vanishing of tautological classes on moduli spaces of curves}, 
preprint, math.AG/0309227.

\bibitem[H]{ha} J. Harer, 
{\em Stability of the homology of the mapping class
groups of orientable surfaces}, Ann. of Math. 121 (1985), 215-249. 

\bibitem[I]{i} E. Ionel, {\em  Topological recursive relations in 
$H^{2g}(\M_{g,n})$}, Invent. Math. 148 (2002), 627-658.

\bibitem[L2]{l1} E. Looijenga, {\em On the tautological ring of $\M_g$}, 
Invent. Math. 121 (1995), 411-419. 

\bibitem[Mo]{mo} S. Morita, {\em Generators for the tautological 
algebra of the moduli space of curves}, Topology 42 (2003), 787-819. 

\bibitem[Mu]{mu} D. Mumford, {\em Towards an enumerative geometry of 
the moduli space of curves} in {\em Arithmetic and geometry II} 
(ed. M. Artin and J. Tate), Progress in Math, vol 36, Birkh\"auser,
 Basel, 1983.

\bibitem[R]{ziv} Z. Ran, {\em Curvilinear enumerative geometry}, Acta
Math. 155 (1985), 81-101.

\end{thebibliography}
\end{document}